\documentclass[10pt,twoside,a4paper,reqno]{amsart}
\usepackage{amscd,amsmath,amsthm,amsfonts,latexsym,amssymb}

\theoremstyle{plain}
\newtheorem{theo+}           {Theorem}
\newtheorem{prop+}           {Proposition}
\newtheorem{coro+}           {Corollary}
\newtheorem{lemm+}           {Lemma}
\newtheorem{conjecture}      {Conjecture}

\theoremstyle{definition}
\newtheorem{defi+}           {Definition}
\newtheorem{problem}         {Problem}
\newtheorem*{ack}            {Acknowledgements}

\theoremstyle{remark}
\newtheorem{rema+}           {Remark}

\newenvironment{theorem}{\begin{theo+}}{\end{theo+}}
\newenvironment{proposition}{\begin{prop+}}{\end{prop+}}
\newenvironment{corollary}{\begin{coro+}}{\end{coro+}}
\newenvironment{lemma}{\begin{lemm+}}{\end{lemm+}}
\newenvironment{remark}{\begin{rema+}}{\end{rema+}}
\newenvironment{definition}{\begin{defi+}}{\end{defi+}}

\newcommand {\bC} {\mathbb {C}}
\newcommand {\bR} {\mathbb {R}}
\newcommand {\bPR} {\mathbb {RP}}
\newcommand {\bN} {\mathbb {N}}
\newcommand {\bx} {\mathbf {x}}
\newcommand {\ze} {\zeta}
\newcommand {\al} {\alpha}
\newcommand {\be} {\beta}
\newcommand {\la} {\lambda}
\newcommand {\si} {\sigma}
\newcommand {\om} {\omega}
\newcommand {\te} {\theta}
\newcommand {\De} {\Delta}
\newcommand {\vf} {\varphi}
\newcommand {\ve} {\varepsilon}
\newcommand {\calL} {\mathcal {L}}
\newcommand {\calH} {\mathcal {H}}
\newcommand {\calC} {\mathcal {C}}
\newcommand {\calD} {\mathcal {D}}
\newcommand {\calR} {\mathcal {R}}
\newcommand {\calA} {\mathcal {A}}
\newcommand {\calB} {\mathcal {B}}
\newcommand {\calZ} {\mathcal {Z}}
\newcommand {\calO} {\mathcal {O}}
\newcommand {\vP} {\varPi}
\newcommand {\E} {\text{End}\,}

\begin{document}

\title[Convexity properties of twisted root maps]{Convexity properties of 
twisted root maps}
\author[J.~Borcea]{Julius Borcea}
\address{Department of Mathematics, Stockholm University, SE-106 91 Stockholm,
   Sweden}
\email{julius@math.su.se}

\begin{abstract}
The strong spectral order induces a natural partial ordering on the manifold 
$\calH_{n}$ of monic hyperbolic polynomials of degree $n$. We prove that 
twisted root maps associated with linear operators acting on $\calH_{n}$ are 
G\aa rding convex on every polynomial pencil and we characterize the class of 
polynomial pencils of logarithmic derivative type by means of the strong 
spectral order. Let ${\calA'}$ be the monoid of linear operators that preserve 
hyperbolicity as well as root sums. We show that any polynomial in 
$\calH_{n}$ is the global minimum of its ${\calA'}$-orbit and we 
conjecture a similar result for complex polynomials.
\end{abstract}

\maketitle

\section*{Introduction}

An important chapter in the theory of distribution of zeros of polynomials 
and transcendental entire functions pertains to the study of linear 
operators that preserve certain prescribed properties 
(cf., e.~g., \cite{CC1,L,RS} and references therein). The following example 
illustrates the viewpoint adopted in this paper. Denote by $\E\vP$ the set 
of linear mappings from the vector space $\vP:=\bC[x]$ to itself and let 
$\vP(\Omega)$ be 
the class of polynomials in $\vP$ whose zeros lie in a fixed set
$\Omega\subseteq \bC$. As noted in \cite{CC1}, the fundamental problem of 
characterizing all operators $T\in \E\vP$ such that 
$T\big(\vP(\Omega)\big)\subseteq \vP(\Omega)$ is open for all but 
trivial choices of $\Omega$. Indeed, this question remains unanswered even in 
the important special cases when $\Omega$ is a line or a half-plane. 
Moreover, in many applications such as stability problems one often 
needs additional information on 
the relative geometry of the zeros of $T(P)$ and $P$ for 
$P\in \vP(\Omega)$ when $T\in \E\vP$ preserves  
$\vP(\Omega)$. For instance, if 
$T=D:=\frac{d}{dx}$ these questions amount to studying the 
geometry of zeros and critical points of complex polynomials, which is 
in itself a vast and intricate subject \cite{RS}. In this case
the Gauss-Lucas theorem implies that the zero set of $T(P)$ is contained in 
the convex hull of the zeros of $P$ and thus 
$T\big(\vP(\Omega)\big)\subseteq \vP(\Omega)$ whenever $\Omega$ is convex. 
However, this result can be substantially refined in various circumstances 
\cite{O,RS}.

In this paper we propose a general setting for studying the relative 
geometry of zeros of polynomials and their distribution under the action of 
various classes of linear operators.
Let $\calC_n$ be the manifold of monic complex polynomials of 
degree $n\ge 1$. For 
$P\in \calC_n$ let $\calZ(P)$ be the unordered $n$-tuple consisting 
of the zeros of $P$, each zero occurring as many times as its multiplicity. 
Hence $\calZ(P)\in \bC^n/\Sigma_n$, where $\Sigma_n$ is the symmetric group on 
$n$ elements. Denote by $\Re \calZ(P)$ (respectively, $\Im \calZ(P)$) the 
unordered $n$-tuple whose components are the real (respectively, imaginary)  
parts of the points in $\calZ(P)$. One says that $P$ 
is {\em hyperbolic} provided that $\Re \calZ(P)=\calZ(P)$, i.~e., if all the 
zeros of $P$ are real. A hyperbolic polynomial with only simple zeros is 
called {\em strictly hyperbolic}. Let $\calH_{n}=\calC_n\cap \vP(\bR)$ be the 
real submanifold of $\calC_n$ 
consisting of hyperbolic polynomials. There is a natural set-theoretic 
identification between $\calC_n$ and $\bC^n/\Sigma_n$ by means of the 
{\em root map}
\begin{equation}\label{rmap}
\begin{split}
\calZ:\calC_n & \longrightarrow \bC^n/\Sigma_n\\
P & \longmapsto \calZ(P)
\end{split}
\end{equation}
whose restriction to $\calH_{n}$ obviously induces a bijection between 
$\calH_{n}$ and $\bR^n/\Sigma_n$. Let $T\in \E\vP$ be an operator such 
that $T(\calC_n)\subseteq \calC_n$. The 
composition $\calZ\circ T$ is called the $T$-{\em twisted root map}. Note 
that if $T$ also acts on $\calH_n$ then the restriction of the $T$-twisted 
root map to $\calH_{n}$ has real components. Given a non-empty set 
$\Omega\subseteq \bC$ we define a multiplicative monoid of linear operators 
by setting
\begin{equation}\label{monon}
\calA_n(\Omega)=\big\{T\in \E\vP\mid 
T\big(\calC_n\cap\vP(\Omega)\big)\subseteq \calC_n\cap\vP(\Omega)\big\}.
\end{equation}
The relevance of twisted root maps in the present context 
is quite clear. Indeed, for degree-preserving linear operators the 
aforementioned questions on the distribution and the relative geometry of 
zeros of polynomials may be summarized as follows.

\begin{problem}\label{pb1}
Describe the properties of $T$-twisted root maps for $T\in \calA_n(\Omega)$, 
where $n$ is a fixed positive integer and $\Omega$ is an appropriate set of 
interest.
\end{problem}

Below we shall mainly focus on the important special case 
of Problem~\ref{pb1} when $\Omega=\bR$.
The following fundamental result from the theory of majorization 
is a key ingredient in our analysis of twisted root maps.

\begin{theorem}\label{spec}
Let $X=(x_{1},x_{2},\ldots,x_{n})^t$ and $Y=(y_{1},y_{2},\ldots,y_{n})^t$ be
two $n$-tuples of vectors in $\bR^{k}$. The following conditions are 
equivalent:
\begin{enumerate}
\item[(i)] For any convex function $f: \bR^k \to \bR$ one
has $\sum_{i=1}^{n}f(x_{i})\le \sum_{i=1}^{n}f(y_{i})$.
\item[(ii)] There exists a doubly stochastic $n\times n$ matrix $A$ such that
$\tilde  X=A\tilde  Y$, where $\tilde X$ and $\tilde Y$ are 
$n\times k$ matrices obtained by some (and then any) ordering of the vectors 
in $X$ and $Y$.
\end{enumerate}
\end{theorem}

If the conditions of Theorem~\ref{spec} are satisfied we say that 
$X$ is {\em majorized} by $Y$ or that $X$ is {\em less than} $Y$ 
{\em in the spectral order}, and write $X\prec Y$. One can easily check that 
$\sum_{i=1}^nx_{i}=\sum_{i=1}^ny_{i}$ if $X\prec Y$. 
Theorem~\ref{spec} is due 
to Schur as well as to Hardy, Littlewood, and P\'olya in the case $k=1$
\cite{HLP}, and to Sherman in the general case \cite{S}. These cases 
are also known as (strong) classical and multivariate majorization, 
respectively. 
Surprisingly, Sherman's theorem was long assumed to be an open problem 
and does not appear in \cite{MO}, which is the definite reference on 
majorization theory (see p.~433 in {\em loc.~cit.}). We refer to \cite{B1} 
for 
a simple new proof of this result. Note that although the spectral order is 
only a 
preordering on $\bR^n$, Birkhoff's theorem \cite[Theorem 2.A.2]{MO} implies 
that it actually induces a partial ordering on $\bR^n/\Sigma_n$. Therefore, 
Theorem~\ref{spec} and the root map in~\eqref{rmap} allow us to 
define a poset structure $(\calH_n,\preccurlyeq)$ by setting 
$P\preccurlyeq Q$ if $P,Q\in \calH_{n}$ and $\calZ(P)\prec \calZ(Q)$.

In \S 1 we establish a general convexity property for 
twisted root maps associated with operators in $\calA_n(\bR)$.
Namely, we show that the restriction of any such map to arbitrary polynomial 
pencils in $\calH_{n}$ is G\aa rding convex 
(Definition~\ref{gard} and Theorem~\ref{GC}). This has several interesting 
consequences for the so-called span (or spread) function and its twisted 
versions (Corollaries~\ref{t-span}-\ref{LD-span}). 

In \S 2 we characterize the class of polynomial pencils of logarithmic 
derivative type contained in $\calH_{n}$ by means of a local minimum property 
with respect to the 
partial ordering $\preccurlyeq$ on $\calH_{n}$ (Theorem~\ref{locm}). 

Let $\calA:=\bigcap_{k=0}^{\infty}\calA_{k}(\bR)$ and denote by $\calA'$ 
the submonoid of $\calA$ whose action on $\vP$ preserves 
the average of the zeros of any polynomial. In \S 3 we show that 
$\calA$ consists of ordinary differential operators of Laguerre-P\'olya type 
(Theorem~\ref{LP}) and that any polynomial in $\calH_{n}$ is the global 
minimum of its $\calA'$-orbit (Theorem~\ref{orbit}). In particular, this 
implies that the action of $\calA$ on $\calH_{n}$ does not decrease the 
span of polynomials (Corollary~\ref{span-A}) and that the polynomial pencils 
characterized in \S 2 satisfy in fact a global minimum
property with respect to the spectral order. 

As we point out in \S 4, a 
natural question that arises from our study is whether one can
describe classical majorization by means of (differential) operators in 
$\calA_n(\bR)$ acting on polynomials in $\calH_{n}$ (Problem~\ref{pb2}). We 
discuss this question as well as possible complex analogs of 
Theorem~\ref{orbit} (Conjecture~\ref{conj1}) and 
extensions of this theorem to the Laguerre-P\'olya 
class of functions.

\begin{ack}
The author would like to thank Petter Br\"and\'en and Boris Shapiro for 
stimulating 
discussions on these and related topics and the anonymous referee for 
useful comments.
\end{ack}

\section{Polynomial pencils and G\aa rding convexity}

A fundamental theorem of G\aa rding asserts that the largest root of a 
multivariate 
homogeneous polynomial which is hyperbolic with respect to a given vector is 
always a convex function \cite{G}. The properties of such polynomials play 
a significant role 
in the theory of partial differential equations, convex analysis and 
matrix theory (see, e.~g., \cite{BGLS}). The following definition is 
motivated by G\aa rding's result.

\begin{definition}\label{gard}
Let $K$ be a convex subset of a vector space. A map 
$f:K\rightarrow \bR^n/\Sigma_n$ given by 
$f(\mathbf{x})=\big(f_1(\mathbf{x}),\ldots,f_n(\mathbf{x})\big)$ is called 
{\em G\aa rding convex} if 
$\mathbf{x}\mapsto \max_{1\le i\le n}f_i(\mathbf{x})$ is a convex 
function on $K$.
\end{definition}

G\aa rding's theorem is a rich source of examples of maps that 
satisfy Definition~\ref{gard} (cf.~\cite{G} and \cite{BGLS}). For instance, 
the restriction 
of the eigenvalue map to the real space of $n\times n$ Hermitian matrices is 
an important such example. 
The main result of this section shows 
that twisted root maps associated to operators in $\calA_n(\bR)$ are 
G\aa rding convex when restricted to certain convex subsets of 
$\calH_n$, as we shall now explain. Recall the notations $\calC_n$, 
$\calH_n$, $\calA_n(\bR)$ from the introduction and denote by $\calR_n$ the 
(real) submanifold of $\calC_n$ consisting of monic real polynomials of 
degree $n$. The inclusion $\calH_n\subseteq \calR_n$ is obviously strict if 
$n\ge 2$, which we assume henceforth.

\begin{definition}\label{pen}
Let $P_1$ and $P_2$ be distinct polynomials in $\calR_n$. The real line 
through $P_1$ and $P_2$, i.~e., the set 
$\calL=\{(1-\la)P_1+\la P_2\mid\la \in \bR\}$, is called a
{\em polynomial pencil} in $\calR_n$. A {\em basis} of a polynomial pencil 
$\calL\subset \calR_n$ is a pair of real 
polynomials $\{P,Q\}$ that satisfy the following conditions:
\begin{equation}\label{bas}
\begin{split}
&P\in \calR_n, \text{ the dominant coefficient of $Q$ equals } n,\\
&\!\deg Q\le n-1, \text{ and } \calL=\{P-\la Q\mid\la \in \bR\}. 
\end{split}
\end{equation}
A polynomial pencil $\calL$ is said to be of 
{\em logarithmic derivative type} or an $LD$-{\em pencil} if there 
exist $Q_1, Q_2, Q_3\in \calL$ such that $Q_3'=Q_1-Q_2$.
\end{definition}

\begin{remark}
As defined above, a polynomial pencil in $\calR_n$ has only 
an affine structure. One can produce real polynomial pencils endowed with
a natural linear structure by using an appropriate projective version 
of Definition~\ref{pen} where a real polynomial pencil is defined instead as 
a real line in projective space $\bPR^n$ identified with the space of all 
homogeneous degree $n$ real polynomials in two real variables considered up to
a constant factor. Such a pencil is called generic if it intersects the
standard discriminant $\calD_{n+1}\subset \bPR^n$ transversally. 
A topological classification of all generic pencils in $\bPR^n$ was 
obtained in \cite{BS1}.
\end{remark}

Clearly, any polynomial pencil has a basis. Moreover, 
if $\{P,Q\}$ and $\{R,S\}$ are two bases of the same polynomial pencil
then $S=Q$ and $R=P-\mu Q$ for some $\mu\in \bR$. Note also that a 
polynomial pencil $\calL\subset \calR_n$ is an $LD$-pencil if and only it has 
a (necessarily 
unique) basis of the form $\{P,P'\}$, which we call the
{\em canonical basis} of the $LD$-pencil $\calL$.
The main result of this section is as follows.

\begin{theorem}\label{GC}
If $\calL$ is an arbitrary polynomial pencil in $\calH_n$ and 
$T\in \calA_n(\bR)$ then the $T$-twisted 
root map $\calZ\circ T|_{\calL}$ is G\aa rding convex.
\end{theorem}

Since we are mainly interested in polynomial pencils contained in 
$\calH_n$ -- which we shall alternatively refer to as {\em hyperbolic} 
(polynomial) 
pencils -- let us first give a complete description of these pencils. 

\begin{theorem}\label{desc}
Let $P$ and $R$ be distinct polynomials in $\calR_n$ and set $Q=P-R$ and 
$\calL=\{P-\la Q\mid\la \in \bR\}$. The following statements are equivalent:
\begin{enumerate}
\item[(i)] $\calL\subset \calH_n$.
\item[(ii)] $\al P+\be R$ is hyperbolic for any $\al, \be\in \bR$ such 
that $\al^2+\be^2\neq 0$.
\item[(iii)] The polynomials $P$ and $R$ are in $\calH_n$ and have weakly 
interlacing zeros.
\item[(iv)] $P\in \calH_n$, $Q$ is hyperbolic, $\deg Q=n-1$, and the zeros of 
$P$ and $Q$ are weakly interlacing. Equivalently, these conditions hold with 
$R$ instead of $P$.
\item[(v)] The polynomial $P+iR$ or indeed $P+iQ$ or $R+iQ$ has all its zeros 
either in the closed upper half-plane or in the closed lower half-plane.
\end{enumerate}
\end{theorem}

\begin{proof}
It is clear that (ii) $\Rightarrow$ (i). Note that if (i) holds then 
$\al P+\be R\in \calH_n$ for any $\al, \be\in \bR$ such 
that $\al^2+\be^2\neq 0$ and $\al +\be\neq 0$. By assumption the polynomial 
$Q$ is not identically zero and in fact $\deg Q\ge 1$, because otherwise 
$P-\la Q\notin \calH_n$ for some $|\la|\gg 0$. Since 
$\{Q-k^{-1}P\}_{k=1}^{\infty}$ is a 
sequence of hyperbolic polynomials which tends to $Q$ uniformly on compact 
sets, Hurwitz' theorem for analytic functions 
\cite[Theorem 1.6.9]{RS} implies that $Q$ is hyperbolic. This proves 
that (i) $\Rightarrow$ (ii). 

The equivalence between (ii) and (iii) is known in the literature 
as Obreschkoff's theorem \cite[Satz 5.2]{O} or the Hermite-Kakeya 
theorem \cite[Theorem 6.3.8]{RS} in the generic case when $P$ and $R$ are 
strictly hyperbolic polynomials with no common zeros. In the general case, 
this equivalence is due to Dedieu \cite[Theorem 4.1]{D}. Actually, the 
arguments used in {\em loc.~cit.~}yield also a proof of 
(i) $\Leftrightarrow$ (iv). Note that the roles of $P$ and $R$ may be 
interchanged and that the condition $\deg Q=n-1$ may alternatively be seen as 
a consequence of (iii) since $P$ and $R$ are distinct.

The well-known Hermite-Biehler theorem \cite[Theorem 6.3.4]{RS} asserts 
that statements (iii) and (v) are equivalent if the words ``weakly'' and 
``closed'' in these statements are replaced by ``strictly'' and ``open'',
respectively.
For the general case, let us set $\tilde{P}=P/S$ and 
$\tilde{R}=R/S$, where $S$ denotes the greatest common monic divisor of 
$P$ and $R$. Note that $\deg \tilde{P}=\deg \tilde{R}\ge 1$ since $P$ and $R$ 
are distinct. If (iii) holds then $\tilde{P}$ and $\tilde{R}$ have strictly 
interlacing zeros (cf.~\cite[Remark 6.3.3]{RS}). By the Hermite-Biehler 
theorem, the polynomial $\tilde{P}+i\tilde{R}$ must have all its zeros 
either in the open upper half-plane or in the open lower half-plane. It 
follows that the polynomial $P+iR=S(\tilde{P}+i\tilde{R})$ has all its zeros 
either in the closed upper half-plane or in the closed lower half-plane, which 
proves (v). Conversely, if the latter holds then $\tilde{P}+i\tilde{R}$ has 
all its zeros either in the open upper half-plane or in the open lower 
half-plane (since any real zero of this polynomial would have to be a common 
zero of $\tilde{P}$ and $\tilde{R}$). Thus $\tilde{P}$ and $\tilde{R}$ have 
strictly interlacing zeros so that $P=S\tilde{P}$ and $R=S\tilde{R}$ have 
weakly interlacing zeros, as stated in (iii). The claims in (v) 
$\Leftrightarrow$ (iv) concerning the polynomials $P$, $Q$, and $P+iQ$ can be 
verified in similar fashion by using the Hermite-Biehler theorem for the pair 
$\{P,Q\}$.
\end{proof}

\begin{proof}[Proof of Theorem~\ref{GC}]
Let $\calL$ be a polynomial pencil 
in $\calH_n$ and $T\in \calA_n(\bR)$. Then either $T|_{\calL}$ is a constant 
map, in which case the conclusion of the theorem holds trivially, or the image 
$T(\calL)$ is again a polynomial pencil in $\calH_n$. Thus, it is enough to 
prove the theorem for $T=Id_{\vP}$, which we assume henceforth. Fix a basis 
$\{P,Q\}$ of $\calL$ as in~\eqref{bas} and denote 
the zeros of $P$ and $Q$ by $x_{1}\le x_{2} \le\ldots \le x_{n}$ and 
$y_{1}\le y_{2} \le \ldots \le y_{n-1}$, respectively. By Theorem~\ref{desc} 
we know that
\begin{equation}\label{interl1}
x_{1}\le y_1\le x_{2} \le\ldots \le x_{n-1}\le y_{n-1}\le x_n.
\end{equation}

Set $R_{\la}=P-\la Q$, $\la\in \bR$, and denote the zeros of $R_{\la}$ by 
$x_{i}(\la)$, $1\le i\le n$, which we label so that 
$x_{i}(0)=x_i$, $1\le i\le n$. Since $R_{\la}-\mu Q\in \calH_n$ for any 
$\mu\in \bR$, Theorem~\ref{desc} again implies that 
\begin{equation}\label{interl2}
x_{1}(\la)\le y_1\le x_{2}(\la) \le\ldots \le x_{n-1}(\la)\le y_{n-1}\le 
x_n(\la)\,\text{ for }\,\la\in \bR.
\end{equation}

{\em Step 1: $P$ and $Q$ have strictly interlacing zeros}. This means that 
$P$ and $Q$ are strictly hyperbolic and have no common zeros. Clearly, any 
common zero of $R_{\la}$ and $Q$ would also have to be a zero of $P$. It 
follows that the interlacing properties in~\eqref{interl1} 
and~\eqref{interl2} are both strict. In 
particular, the polynomial $R_{\la}$ is strictly hyperbolic for any 
$\la\in \bR$. We may therefore differentiate the identities 
$R_{\la}(x_i(\la))=0$, $1\le i\le n$, with respect to $\la$ to get
\begin{equation}\label{deriv1}
x_{i}'(\la)=\frac{Q(x_i(\la))}{R_{\la}'(x_i(\la))}\,\text{ for }\,
\la\in \bR\,\text{ and }\,1\le i\le n,
\end{equation}
where $R_{\la}'(x)=\frac{\partial}{\partial x}R_{\la}(x)$. Note 
that~\eqref{deriv1} readily 
implies that $x_{i}'(\la)>0$ for all $\la\in \bR$ and $1\le i\le n$ since 
both $Q(x_i(\la))$ and $R_{\la}'(x_i(\la))$ have constant signs while
$$\frac{Q(x_i(0))}{P'(x_i(0))}=\frac{Q(x_i)}{P'(x_i)}=
n\prod_{j=1}^{n-1}(x_i-y_j)
\prod_{\substack{j=1\\ j\neq i}}^{n}(x_i-x_j)^{-1}>0$$
because of the (strict) inequalities in~\eqref{interl1}. Thus, all 
the zeros of the polynomial $R_{\la}$ are increasing functions of $\la$. 
By differentiating~\eqref{deriv1} with respect to $\la$ we obtain
\begin{equation}\label{deriv2}
\begin{split}
x_{i}''(\la) & = x_{i}'(\la)\left[\frac{Q'R_{\la}'-QR_{\la}''}
{R_{\la}'^2}\right]\!(x_i(\la))+\left[\frac{QQ'}{R_{\la}'^2}\right]\!
(x_i(\la))= 2x_{i}'(\la)\frac{Q'(x_i(\la))}{R_{\la}'(x_i(\la))}\\
& \phantom{= x} -\big(x_{i}'(\la)\big)^2
\frac{R_{\la}''(x_i(\la))}{R_{\la}'(x_i(\la))}
= \big(x_{i}'(\la)\big)^2\left[\frac{2Q'}{Q}-
\frac{R_{\la}''}{R_{\la}'}\right]\!(x_i(\la))
\end{split}
\end{equation}
for $\la\in \bR$ and $1\le i\le n$. The special case when $i=n$ 
in~\eqref{deriv2} yields
\begin{equation*}
\begin{split}
x_{n}''(\la) & = \big(x_{n}'(\la)\big)^2\left[\sum_{j=1}^{n-1}
\frac{2}{x_n(\la)-y_j}-\sum_{j=1}^{n-1}\frac{2}{x_n(\la)-x_j(\la)}\right]\\
& = 2\big(x_{n}'(\la)\big)^2\sum_{j=1}^{n-1}\frac{y_j-x_j(\la)}
{(x_n(\la)-y_j)(x_n(\la)-x_j(\la))}>0
\end{split}
\end{equation*}
by~\eqref{interl2}. This implies that 
$\bR\ni \la\mapsto \max \calZ(R_{\la})=x_n(\la)$ is a convex function. Thus 
$\calZ|_{\calL}$ is a G\aa rding convex map, which 
proves the theorem in the generic case.

{\em Step 2: The general case.} Let $S$ denote the greatest common monic 
divisor of $P$ and $Q$. Step 1 shows that the theorem is true if $S\equiv 1$ 
and so we may assume that $\deg S\ge 1$. Set $\tilde{P}=P/S$, 
$\tilde{Q}=Q/S$, and $\tilde{R}_{\la}=\tilde{P}-\la \tilde{Q}$, so that 
$R_{\la}=S\tilde{R}_{\la}$ and thus
$\max \calZ(R_{\la})=\max \big(\al, \max \calZ(\tilde{R}_{\la})\big)$ for 
$\la \in \bR$, where 
$\al$ is the largest zero of $S$. If $\tilde{Q}$ is a constant polynomial 
(which by \eqref{bas} would actually mean that $\tilde{Q}\equiv n$) then 
$\max \calZ(\tilde{R}_{\la})$ is obviously a linear function of $\la$. 
Otherwise $\tilde{P}$ and $\tilde{Q}$ must have strictly interlacing zeros 
and so by step 1 the function $\bR\ni \la\mapsto \max \calZ(\tilde{R}_{\la})$ 
is convex.
In either case it follows that $\bR\ni \la\mapsto \max \calZ(R_{\la})$ is a 
convex function so that $\calZ|_{\calL}$ is a G\aa rding convex map, which 
completes the proof.
\end{proof}

\begin{remark}\label{rem1}
Theorem~\ref{GC} generalizes the result announced in Theorem 1.7 of \cite{BS2},
where the G\aa rding convexity property was stated only for $LD$-pencils in 
$\calH_n$.  
\end{remark}

Recall that the span (or spread) of a polynomial $P\in \calH_n$ is 
the length of the smallest interval that contains all its zeros, 
i.~e., $\De(P)=\max \calZ(P)-\min \calZ(P)$. A review of the literature on 
the span of hyperbolic polynomials and related questions may be found in 
\cite[Ch.~6]{RS}. Given an operator $T\in \calA_n(\bR)$ we define the 
$T$-{\em twisted span function} on $\calH_n$ to be the composite map 
$\De\circ T$. From Theorem~\ref{GC} we deduce the following 
properties for twisted span functions:

\begin{corollary}\label{t-span}
If $\calL$ is an arbitrary polynomial pencil in $\calH_n$ and 
$T\in \calA_n(\bR)$ then the $T$-twisted span function 
$\De\circ T|_{\calL}$ is convex and has a global minimum. 
\end{corollary}

\begin{proof}
As in the proof of Theorem~\ref{GC}, it is enough to consider the case 
$T=Id_{\vP}$. Let $\calL$ be a polynomial pencil in $\calH_n$ with a basis 
$\{P,Q\}$ as in \eqref{bas}. Set $\tilde{P}(x)=(-1)^nP(-x)$ and 
$\tilde{Q}(x)=(-1)^{n-1}Q(-x)$, so that $\calZ(P-\la Q)=
-\calZ(\tilde{P}+\la \tilde{Q})$ for all $\la\in \bR$. By Theorem~\ref{GC} the 
function $\bR\ni \la\mapsto \max\calZ(\tilde{P}+\la \tilde{Q})$ is convex 
and thus $\bR\ni \la\mapsto \min\calZ(P-\la Q)=-\max\calZ(\tilde{P}+
\la \tilde{Q})$ is a concave function. Hence $\bR\ni \la\mapsto \De(P-\la Q)$ 
is a convex function and therefore also Lipschitz continuous on any compact 
interval. Since $\De(P-\la Q)\rightarrow \infty$ as $|\la|\rightarrow \infty$ 
it follows that $\De|_{\calL}$ has a global minimum, as required.
\end{proof}

Corollary~\ref{t-span} can be further refined in the case of hyperbolic 
$LD$-pencils:

\begin{corollary}\label{LD-span}
If $\calL$ is an $LD$-pencil in $\calH_n$ with canonical basis $\{P,P'\}$ 
then the span function $\De|_{\calL}$ is convex and has a global minimum at 
$P$. In particular, for any $P\in \calH_n$ and $\la\in \bR$ one has 
$\De(P)\le \De(P-\la P')$. More generally, if $\la_1, \la_2\in \bR$ are such 
that $\la_1 \la_2\ge 0$ and $|\la_1|\le |\la_2|$ then 
$\De(P-\la_1 P')\le \De(P-\la_2 P')$. 
\end{corollary}

\begin{proof}
It is clearly enough to prove only the last assertion of the corollary. Let 
$P$ be 
a strictly hyperbolic polynomial in $\calH_n$ with zeros $x_1<\ldots<x_n$. 
Denote by $x_i(\la)$, $1\le i\le n$, the zeros of the strictly hyperbolic 
polynomial $P-\la P'$, $\la\in \bR$, and assume that these are labeled so that 
$x_i(0)=x_i$, $1\le i\le n$. The arguments in the proof of Theorem~\ref{GC} 
show that $x_i(\la)<x_{i+1}(\la)$ for all $\la\in \bR$ and $1\le i\le n-1$. 
Moreover, by using \eqref{deriv1} and \eqref{deriv2} with $Q=P'$ we see that 
$\bR\ni \la\mapsto x_n(\la)-\la$ is a strictly convex function with a global 
minimum at $\la=0$ while $\bR\ni \la\mapsto x_1(\la)-\la$ is a (strictly) 
concave function with a global maximum at $\la=0$, which proves the corollary 
in the generic case when $P$ has simple zeros. If $P$ has multiple zeros then 
we consider the strictly hyperbolic polynomial $P_{\ve}$ with zeros 
$x_k+k\ve$, $1\le k\le n$, where $\ve>0$. Let $\la_1, \la_2\in \bR$ be as 
in the corollary. The desired conclusion follows by letting 
$\ve\rightarrow 0$ in the inequality 
$\De(P_{\ve}-\la_1 P_{\ve}')\le \De(P_{\ve}-\la_2 P_{\ve}')$, which holds 
thanks to the first part of the proof since $P_{\ve}$ has only simple zeros.
\end{proof}

The study of the geometrical structure of $\calH_n$ was initiated by 
Arnold in \cite{A}. Subsequently, the convex subsets of $\calH_n$ were 
characterized in \cite{D}. In view of the above results one may ask whether 
twisted root maps are G\aa rding convex when restricted to arbitrary convex 
subsets of $\calH_n$. This is definitely not true, as one can see by 
considering for instance the subset
$\{(1-\al)P_1+\al P_2\mid \al\in [0,1]\}$ of $\calH_2$, where 
$P_1(x)=x^2-1$ and $P_2(x)=x^2-2$. However, we can show that the following 
analog of Theorem~\ref{GC} holds for arbitrary segments of hyperbolic 
polynomials:

\begin{corollary}\label{segm}
Let $P_1$ and $P_2$ be distinct polynomials in $\calH_n$ such that the segment 
$[P_1,P_2]:=\{(1-\al)P_1+\al P_2\mid \al\in [0,1]\}$ is contained in 
$\calH_n$. There 
exists a (non-unique) polynomial $P_3\in \calH_n$ such that for any 
$T\in \calA_n(\bR)$ the $T$-twisted root maps $\calZ\circ T|_{[P_1,P_3]}$ and 
$\calZ\circ T|_{[P_2,P_3]}$ are G\aa rding convex.
\end{corollary}

\begin{proof}
Let $x_1\le\ldots\le x_n$ and $y_1\le\ldots \le y_n$ be the zeros of $P_1$ 
and $P_2$, respectively. According to \cite[Theorem 2.1]{D}, the segment 
$[P_1,P_2]$ is contained in $\calH_n$ if and only if $\max(x_i,y_i)\le 
\min(x_{i+1},y_{i+1})$ for $1\le i\le n-1$. Let $z_i$, $1\le i\le n$, be such 
that $z_n\ge \max(x_n,y_n)$ and 
$\max(x_i,y_i)\le z_i\le\min(x_{i+1},y_{i+1})$, $1\le i\le n-1$. 
Denote by $\calL_k$, $k=1,2$, the real line through $P_k$ and $P_3$, $k=1,2$,
where $P_3\in \calH_n$ is such that $\calZ(P_3)=(z_1,\ldots,z_n)$. Since both
pairs of polynomials $\{P_1,P_3\}$ and $\{P_2,P_3\}$ have weakly interlacing 
zeros Theorem~\ref{desc} implies that both $\calL_1$ and 
$\calL_2$ are polynomial pencils contained in $\calH_n$. The result follows 
readily from Theorem~\ref{GC}.
\end{proof}

\section{A characterization of hyperbolic $LD$-pencils}

Hyperbolic pencils of logarithmic derivative type are particularly interesting
for at least two reasons. On the one hand, they are obviously related to 
the action of linear differential operators on the manifold $\calH_n$. On 
the other hand, these pencils have interesting connections 
with classical majorization via the partial ordering on $\calH_n$ defined in 
the introduction. Indeed, the main result below states that the 
class of hyperbolic $LD$-pencils is actually characterized by a local minimum 
property with respect to the spectral order. This hints at possibly even 
deeper connections between hyperbolic polynomials, classical majorization, 
and differential operators, which we shall further investigate in the next 
sections. 

\begin{definition}\label{shift}
Let $\calL$ be a polynomial pencil in $\calR_n$ with a basis $\{P,Q\}$ as 
in \eqref{bas}. The set 
$\calL_s(P,Q):=\{P(x+\la)-\la Q(x+\la)\mid\la \in \bR\}$ is called the 
$\{P,Q\}$-{\em shift} of $\calL$.
\end{definition}

Note that all polynomials in the $\{P,Q\}$-shift of $\calL$  have the same 
zero sum as $P$.

\begin{theorem}\label{locm}
A polynomial pencil $\calL\subset \calH_n$ is an $LD$-pencil if and only if 
there is a shift $\calL_s(P,Q)$ of $\calL$ such that the root map 
$\calZ|_{\calL_s(P,Q)}$ has a local minimum with respect to the spectral 
order, i.~e., there exist $\la_0\in \bR$ and $\ve>0$ such that 
$$P(x+\la_0)-\la_0 Q(x+\la_0)\preccurlyeq P(x+\la)-\la Q(x+\la)\,\text{ for }\,
\la\in (\la_0-\ve,\la_0+\ve),$$
where $\preccurlyeq$ denotes the partial ordering on $\calH_n$.
\end{theorem}

The proof of the sufficiency part of Theorem~\ref{locm} relies on the 
following lemma.

\begin{lemma}\label{suff}
Let $\calL$ be a polynomial pencil in $\calH_n$ with a basis $\{P,Q\}$ as 
in \eqref{bas} and assume that there exists $\ve>0$ such that 
\begin{equation}\label{min}
P(x)\preccurlyeq P(x+\mu)-\mu Q(x+\mu)\,\text{ for }\,\mu\in (-\ve,\ve).
\end{equation}
Then $Q=P'$, so that $\calL$ is an $LD$-pencil with canonical basis $\{P,P'\}$.
\end{lemma}

\begin{proof}
Let $S$ denote the greatest common monic divisor of $P$ and $Q$. Then either 
$S\equiv 1$ or $\deg S\ge 1$, in which case we denote 
by $w_1,\ldots,w_k$ the distinct zeros of $S$ with multiplicities
$s_1,\ldots,s_k$, respectively. Hence we may write
$$S(x)=\prod_{i=1}^{k}(x-w_i)^{s_i},\,\,\tilde{P}(x):=\frac{P(x)}{S(x)}
=\prod_{i=1}^{d}(x-x_i),
\,\,\tilde{Q}(x):=\frac{Q(x)}{S(x)}=n\prod_{i=1}^{d-1}(x-y_i),$$
where $d\ge 1$, $\sum_{i=1}^{k}s_i=n-d$, $x_1<y_1<x_2<\ldots<x_{d-1}<
y_{d-1}<x_d$ if $d\ge 2$, with the usual understanding that empty products
are equal to one while empty sums equal zero. For $\mu\in \bR$ let 
$\tilde{R}_{\mu}=\tilde{P}-\la\tilde{Q}$ and
$R_{\mu}=P-\mu Q=S\tilde{R}_{\mu}$. Denote by $x_i(\mu)$, 
$1\le i\le d$, the zeros of $\tilde{R}_{\mu}$, which we 
label so that $x_{i}(0)=x_i$ for $1\le i\le d$. Since $\tilde{R}_{\mu}$ and 
$\tilde{Q}$ have strictly interlacing zeros if $d\ge 2$, all these zeros are 
simple. We may therefore use a computation similar to \eqref{deriv1} to get
\begin{equation}\label{deriv3}
x_{i}'(\mu)=\frac{\tilde{Q}(x_i(\mu))}{\tilde{R}_{\mu}'(x_i(\mu))}
\,\text{ for }\,\mu\in \bR\,\text{ and }\,1\le i\le d.
\end{equation}
Choose $c\in \bR$ such that 
$\ze+c\neq 0$ whenever $P(\ze)=0$ and consider the sequence of convex 
functions $\{f_m\}_{m=1}^{\infty}$ given by $f_m(x)=(x+c)^m$. Note that the 
complete list of zeros of $R_{\mu}$ consists of $w_j-\mu$, $1\le j\le k$, 
with multiplicities $s_1,\ldots,s_k$, respectively, and $x_i(\mu)-\mu$, 
$1\le i\le d$, and that condition \eqref{min} reads 
$R_0\preccurlyeq R_{\mu}$ for $|\mu|<\ve$. By Theorem~\ref{spec} this implies
that the differentiable function
$$\bR\ni \mu\mapsto \sum_{i=1}^{d}f_{m}(x_{i}(\mu)-\mu)+
\sum_{j=1}^{k}s_{j}f_{m}(w_j-\mu)$$
has a local minimum at $\mu=0$ for any fixed $m\in \bN$. Differentiation with 
respect to
\noindent
$\mu$ and formula \eqref{deriv3} then yield the identities
\begin{equation}\label{deriv4}
\sum_{i=1}^{d}\left(\frac{\tilde{Q}(x_i)}{\tilde{P}'(x_i)}-1\right)
(x_{i}+c)^{m-1}
-\sum_{j=1}^{k}s_{j}(w_{j}+c)^{m-1}=0\,\text{ for }m\in \bN.
\end{equation}
From \eqref{deriv4} and the choice of $c$ we deduce that the following 
relations must hold:
\begin{equation*}
\begin{split}
& k\le d;\,\,\tilde{P}(w_j)=0\,\text{ for }1\le j\le k;
\,\,P(x)=\prod_{i=1}^{d}(x-x_i)^{n_i},\text{ where }n_i=s_i+1\\
& \text{if }S(x_i)=0\text{ and }n_i=1\text{ otherwise};\,\text{ and }
\frac{\tilde{Q}(x_i)}{\tilde{P}'(x_i)}=n_i\,\text{ for }1\le i\le d.
\end{split}
\end{equation*}
Using these relations and a partial fractional decomposition we obtain
$$\frac{Q(x)}{P(x)}=\frac{S(x)\tilde{Q}(x)}{S(x)\tilde{P}(x)}
=\frac{\tilde{Q}(x)}{\tilde{P}(x)}=\sum_{i=1}^{d}\frac{\tilde{Q}(x_i)}
{\tilde{P}'(x_i)}\frac{1}{x-x_i}=\sum_{i=1}^{d}\frac{n_i}{x-x_i}
=\frac{P'(x)}{P(x)}$$
for all $x\neq x_i$, $1\le i\le d$. It follows that $Q=P'$, as required.
\end{proof}

The proof of the necessity part of Theorem~\ref{locm} is based on a
criterion for classical majorization due to Hardy, 
Littlewood, and P\'olya \cite{HLP}. It should be mentioned that there are no 
known analogs of this criterion for multivariate majorization.

\begin{theorem}\label{crit}
Let $X=(x_{1}\le x_{2} \le\ldots \le x_{n})$ and
$Y=(y_{1}\le y_{2} \le \ldots \le y_{n})$ be two $n$-tuples of 
real numbers. Then $X\prec Y$ if and only if the $x_{i}$'s and the 
$y_{i}$'s 
satisfy the following conditions:
$$\sum_{i=1}^{n}x_{i}=\sum_{i=1}^{n}y_{i}\,\text{ and }\,
\sum_{i=0}^{k}x_{n-i}\le \sum_{i=0}^{k}y_{n-i}\,\text{ for }\,
0\le k\le n-2.$$
\end{theorem}

\begin{lemma}\label{nec}
If $P\in \calH_n$ then there exists $\ve>0$ such that for all real $\la$ with 
$|\la|<\ve$ one has $P(x)\preccurlyeq P(x+\la)-\la P'(x+\la)$.
\end{lemma}

\begin{proof}
We show first that for any $P\in \calH_n$ there exists $\ve_1=\ve_1(P)>0$ 
such that if $\la\in [0,\ve_1)$ then $P(x)\preccurlyeq P(x+\la)-\la P'(x+\la)$.
Let $x_1<\ldots<x_d$ denote the distinct zeros of $P$ with multiplicities 
$n_1,\ldots,n_d$, respectively. Set 
$S(x)=\prod_{i=1}^{d}(x-x_i)^{n_i-1}$, $\tilde{P}=P/S$, $Q=P'/S$, and 
$\tilde{R}_{\la}=\tilde{P}-\la Q$, $\la\in \bR$. Clearly, $\tilde{R}_{\la}$ 
is a strictly hyperbolic polynomial for all $\la\in \bR$. Denote its 
zeros by $x_1(\la),\ldots,x_d(\la)$ and label these so that $x_i(0)=x_i$, 
$1\le i\le d$. An argument similar to the one used in the proof of 
Theorem~\ref{GC} shows that all these zeros are increasing functions of 
$\la$. Moreover, by analogy with \eqref{deriv1} and \eqref{deriv2} and 
some straightforward computations we obtain
\begin{equation}\label{approx}
\begin{split}
& x_{i}'(0)=\frac{Q(x_i(0))}{\tilde{P}'(x_i(0))}=
\frac{Q(x_i)}{\tilde{P}'(x_i)}=n_i\text{ and }
x_{i}''(0)=\big(x_{i}'(0)\big)^2\left[\frac{2Q'}{Q}-
\frac{\tilde{R}_{0}''}{\tilde{R}_{0}'}\right]\!(x_{i}(0))\\
& \!=n_{i}^{2}\left[\frac{2Q'(x_i)}{Q(x_i)}-
\frac{\tilde{P}''(x_i)}{\tilde{P}'(x_i)}\right]
=2\sum_{\substack{j=1\\ j\neq i}}^{d}\frac{n_{i}n_{j}}{x_i-x_j}
\,\text{ for }1\le i\le d;\,\,\sum_{i=1}^{d}x_{i}''(0)=0;\\
&\text{and }\sum_{i=1}^{k}x_{i}''(0)=2\sum_{i=1}^{k}\sum_{j=k+1}^{d}
\frac{n_{i}n_{j}}{x_i-x_j}<0\text{ if }k\le d-1\text{ since }
x_1<\ldots<x_d.
\end{split}
\end{equation}
Assume for now that $\la\ge 0$ and set $R_{\la}=P-\la P'=S\tilde{R}_{\la}$. 
Let $z_m(\la)$, $1\le m\le n$, be the zeros of $R_{\la}$, which we label as 
follows. Given $m\in \{1,2,\ldots,n\}$ there is a unique 
$i=i(m)\in \{1,\ldots,d\}$ 
such that $\sum_{j=0}^{i-1}n_j<m\le \sum_{j=0}^{i}n_j$, where $n_0:=0$. Then 
we set $z_m(\la)=x_i$, $\la\ge 0$, if $n_i\ge 2$ and $m<\sum_{j=0}^{i}n_j$, 
and we let $z_m(\la)=x_i(\la)$, $\la\ge 0$, otherwise. Note 
that with this labeling we have $z_1(\la)\le z_2(\la)\le \ldots\le z_n(\la)$
for any $\la\ge 0$. Furthermore, if $1\le m\le n-1$ then for all small 
$\la\ge 0$ we get
\begin{equation}\label{approx1}
\begin{split}
& \sum_{j=1}^{m}\big(z_j(\la)-\la\big)-\sum_{j=1}^{m}z_j(0)=
-\left(\!m-\sum_{j=0}^{i(m)-1}\!\!\!n_j\!\right)\!\la+\calO(\la^2)
\,\text{ if }\, m<\sum_{j=0}^{i(m)}n_j,\\
& \sum_{j=1}^{m}\big(z_j(\la)-\la\big)-\sum_{j=1}^{m}z_j(0)=\frac{1}{2}\!
\left(\sum_{j=1}^{i(m)}x_{j}''(0)\right)\!\la^2+\calO(\la^3)\,\text{ if }\, 
m=\sum_{j=0}^{i(m)}n_j.
\end{split}
\end{equation}
From \eqref{approx} and \eqref{approx1} we see that there exists 
$\ve_1=\ve_1(P)>0$ such that if $\la\in [0,\ve_1)$ and $1\le m\le n-1$ then
$\sum_{j=1}^{m}(z_j(\la)-\la)\le \sum_{j=1}^{m}z_j(0)$, which is the same as 
$$\sum_{i=1}^{k}z_{n-i}(0)\le \sum_{i=1}^{k}\big(z_{n-i}(\la)-\la\big)
\text{ for }\la\in [0,\ve_1)\text{ and } 0\le k\le n-2$$ 
since $\sum_{j=1}^{n}(z_j(\la)-\la)=\sum_{j=1}^{n}z_j(0)$ whenever $\la\ge 0$. 
By Theorem~\ref{crit} this means that 
$P(x)\preccurlyeq P(x+\la)-\la P'(x+\la)$ for $\la\in [0,\ve_1)$, as required.
In order to complete the proof let $P_1(x)=(-1)^{n}P(-x)\in \calH_n$. 
The above arguments applied to $P_1$ show that there exists some 
$\ve_2=\ve_2(P)>0$ such 
that $P_1(x)\preccurlyeq P_1(x+\mu)-\mu P_{1}'(x+\mu)$ for all 
$\mu\in [0,\ve_2)$. Since $\calZ\big(P_1(x+\mu)-\mu P_{1}'(x+\mu)\big)=
-\calZ\big(P(x+\la)-\la P'(x+\la)\big)$, where $\la=-\mu$, it follows that 
$P(x)\preccurlyeq P(x+\la)-\la P'(x+\la)$ for any real $\la$ with 
$|\la|<\ve:=\min (\ve_1,\ve_2)$. This finishes the proof of the lemma.
\end{proof}

\begin{proof}[Proof of Theorem~\ref{locm}] 
Let $\calL$ be a polynomial pencil in $\calH_n$ such that there exists a shift 
$\calL_s(P,Q)$ of $\calL$ that satisfies the local minimum property stated in 
the theorem with $\la_0\in \bR$ and $\ve>0$. Set 
$\hat{P}(x)=P(x+\la_0)-\la_0 Q(x+\la_0)$, $\hat{Q}(x)=Q(x+\la_0)$, and 
$\mu=\la-\la_0$. Clearly, the local minimum condition translates into 
$$\hat{P}(x)\preccurlyeq \hat{P}(x+\mu)-\mu \hat{Q}(x+\mu)\,\text{ for }\, 
\mu\in (-\ve,\ve).$$ 
Applying Lemma~\ref{suff} to the polynomial pencil 
$\hat{\calL}:=\{\hat{P}-\mu \hat{Q}\mid \mu\in \bR\}$ we get 
$\hat{Q}=\hat{P}'$. Hence $Q=P'-\la_0 Q'$, so that   
$(P-\la_0 Q)'=(P-\la_0 Q)-(P-(\la_{0}+1)Q)$ and thus $\calL$ is an 
$LD$-pencil by Definition~\ref{pen}. Conversely, if $\calL$ is an 
$LD$-pencil in $\calH_n$ with canonical basis $\{P,P'\}$ then Lemma~\ref{nec} 
shows that the shift $\calL_s(P,P')$ satisfies the local minimum property 
stated in the theorem.
\end{proof}

\begin{remark}\label{gen-LD}
The necessity part of Theorem~\ref{locm} may also be seen as 
a corollary of Theorem~\ref{orbit} below, where it is shown that 
$LD$-pencils in $\calH_n$ satisfy in fact a global minimum property with 
respect to the spectral order. Note also that unlike the 
property stated in Theorem~\ref{locm}, the minimum property for 
(twisted) span functions obtained in Corollary~\ref{t-span} is not specific 
for the class of hyperbolic $LD$-pencils.
\end{remark}

\section{Spectral order and differential operators \\ of Laguerre-P\'olya type}

The monoid $\calA_n(\bR)=\big\{T\in \E\vP\mid T\big(\calH_n\big)\subseteq 
\calH_n\big\}$ was previously defined only for $n\ge 1$. Let us 
extend this notation to $n=0$ by putting $\calH_0=\{1\}\subset \vP$ and 
$\calA_0(\bR)=\big\{T\in \E\vP\mid T(\calH_0)=\calH_0\}$. Given a 
non-constant polynomial 
$P\in \vP$ we denote by $\si(P)$ the sum of the zeros of $P$. Set
$$\calA=\bigcap_{n=0}^{\infty}\calA_n(\bR)\,\text{ and }\,
\calA'=\big\{T\in \calA\mid \si\big(T(P)\big)=\si(P)\text{ if } P\in \vP,
\,\deg P\ge 1\big\}.$$
Thus $\calA$ is the largest monoid of linear operators that act on each 
$\calH_n$ for $n\ge 0$ while $\calA'$ is the largest submonoid of 
$\calA$ consisting of operators whose action on $\vP$ preserves the average 
of the zeros of any non-constant polynomial. As we already saw in \S 2, 
hyperbolic $LD$-pencils may be described by means of a local minimum 
property that involves the spectral order on $\bR$. More generally, the study 
of the relative location of the zero sets $\calZ\big(T(P)\big)$ and $\calZ(P)$ 
for $T\in \calA'$ and $P\in \calH_n$ reveals some interesting connections 
between the action of hyperbolicity-preserving linear operators on the 
manifold $\calH_n$ and classical majorization. Indeed, the main result of this 
section shows that if $n\ge 1$ then any polynomial in $\calH_n$ is the global 
minimum of its $\calA'$-orbit with respect to the partial ordering on 
$\calH_n$:

\begin{theorem}\label{orbit}
If $n\ge 1$ and $P\in \calH_n$ then $P\preccurlyeq T(P)$ for any $T\in \calA'$.
\end{theorem}

Before embarking on the proof let us point out that 
\cite[Theorem~1]{CPP} and the Hermite-Poulain theorem yield actually a 
complete description of the structure of the monoids $\calA$ and $\calA'$:

\begin{theorem}\label{LP}
An operator $T\in \E\vP$ belongs to $\calA$ if and only if $T=\vf(D)$, where 
$D=d/dx$ and $\vf$ is a real entire function in the Laguerre-P\'olya class of 
the form
$$\vf(x)=e^{-a^{2}x^{2}+bx}\prod_{k=1}^{\infty}(1-\al_{k}x)e^{\al_{k}x}$$
with $a,b,\al_{k}\in \bR$ and $\sum_{k=1}^{\infty}\al_{k}^2<\infty$. In 
particular, $\calA$ is a commutative monoid.
\end{theorem}

\begin{corollary}\label{LP'}
The monoid $\calA'$ consists of linear operators of the form $\vf(D)$, 
where $D=d/dx$ and $\vf$ is a real entire function in the Laguerre-P\'olya 
class given by 
$$\vf(x)=e^{-a^{2}x^{2}}\prod_{k=1}^{\infty}(1-\al_{k}x)e^{\al_{k}x}$$
with $a,\al_{k}\in \bR$ and $\sum_{k=1}^{\infty}\al_{k}^2<\infty$. Thus 
$\calA=\calA'\times \left\langle e^{bD}\mid b\in \bR\right\rangle$.
\end{corollary}

In view of Theorem~\ref{LP} it seems reasonable to adopt the following 
terminology: an operator $T\in \E\vP$ 
is said to be a {\em differential operator of Laguerre-P\'olya type} if 
$T=\vf(D)$, where $D=d/dx$ and $\vf$ is a real entire function in the 
Laguerre-P\'olya class. Such operators were studied in e.~g.~\cite{CC2} in 
connection with various generalizations of the P\'olya-Wiman 
conjecture. Since $\calA'$ contains only differential 
operators of Laguerre-P\'olya type, it is enough to check that 
Theorem~\ref{orbit} is true for the ``building blocks'' of these operators, 
that is, differential operators of the form $(1-\la D)e^{\la D}$ or 
$e^{-\la^{2} D}$ with $\la\in \bR$. To do this we need the following lemma.

\begin{lemma}\label{ind}
Let $n\ge 2$, $\bx=(x_1<\ldots<x_n)\in \bR^n$, and 
$P(x)=\prod_{i=1}^{n}(x-x_i)\in \calH_n$. For $\la\in \bR$ denote by 
$\ze_i=\ze_{i}(\la;\bx)$, $1\le i\le n$, the zeros of the strictly hyperbolic 
polynomial $P-\la P'$. If these are labeled so that 
$\ze_{i}(0;\bx)=x_i$, $1\le i\le n$, then
$$\big(\ze_{i}(\la;\bx)-x_j\big)^{2}
\frac{\partial \ze_i}{\partial x_j}(\la;\bx)
=\la^{2}\frac{\partial \ze_i}{\partial \la}(\la;\bx)>0\text{ for }
\la\neq 0\text{ and }1\le i,j\le n.$$
In particular, for any fixed values $x_1<\ldots<x_{n-1}$ and $\la\neq 0$ each 
of the functions $(x_{n-1},\infty)\ni x_n\mapsto \ze_{i}(\la;\bx)$, 
$1\le i\le n$, is increasing.
\end{lemma}

\begin{proof}
Let $1\le i,j\le n$, $\bx\in \bR^n$, and set $P(x)=(x-x_j)Q(x)$. For 
$\la\in \bR$ we get
\begin{equation*}
\begin{split}
& P(x)-\la P'(x)=(x-x_j)\left[Q(x)-\la Q'(x)\right]-\la Q(x),\\
& P'(x)-\la P''(x)
=Q(x)-2\la Q'(x)+(x-x_j)\left[Q'(x)-\la Q''(x)\right],\\
& P\big(\ze_{i}(\la;\bx)\big)=\la P'\big(\ze_{i}(\la;\bx)\big)
=\big(\ze_{i}(\la;\bx)-x_j\big)Q\big(\ze_{i}(\la;\bx)\big),\\
& \text{and }\!\left[P'\big(\ze_{i}(\la;\bx)\big)-
\la P''\big(\ze_{i}(\la;\bx)\big)\right]
\frac{\partial \ze_i}{\partial \la}(\la;\bx)=P'\big(\ze_{i}(\la;\bx)\big).
\end{split}
\end{equation*}
The arguments in the the proof of Theorem~\ref{GC} show that
if $\bx$ is fixed then for any $\la$ one has
$\frac{\partial \ze_i}{\partial \la}(\la;\bx)>0$ and $\ze_{i}(\la;\bx)\neq 
x_k$, $1\le k\le n$, if $\la\neq 0$. By 
differentiating the identity
$\big(\ze_{i}(\la;\bx)-x_j\big)
\left[Q\big(\ze_{i}(\la;\bx)\big)-\la Q'\big(\ze_{i}(\la;\bx)\big)\right]
=\la Q\big(\ze_{i}(\la;\bx)\big)$ with respect to $x_j$ 
and using the relations listed above we arrive at the desired conclusion.
\end{proof}

\begin{proof}[Proof of Theorem~\ref{orbit}]
The theorem holds trivially for $n=1$ since $T|_{\calH_1}=Id_{\calH_1}$ if 
$T\in \calA'$ by Corollary~\ref{LP'}. Hence we may assume that $T\in \calA'$ 
and $P\in \calH_n$ with $n\ge 2$. 

{\em Step 1: $P$ is strictly hyperbolic and $T=(1-\la D)e^{\la D}$ for some 
$\la\in \bR$.} Using the notations of Lemma~\ref{ind} we denote the zeros of 
$P$ by $x_1<\ldots<x_n$ and those of the (strictly hyperbolic) polynomial 
$P_{\la}:=P-\la P'$ by $\ze_i=\ze_{i}(\la;\bx)$, $1\le i\le n$, where 
$\bx=(x_1<\ldots<x_n)\in \bR^n$. We further assume that the latter are 
labeled so that $\ze_{i}(0;\bx)=x_i$, $1\le i\le n$. As in the first step   
of the proof of Theorem~\ref{GC} we see that this labeling of the zeros yields 
$\ze_{1}(\la;\bx)<\ldots<\ze_{n}(\la;\bx)$ for all $\la\in \bR$. By 
Theorem~\ref{crit} the relation $P\preccurlyeq T(P)$ is equivalent to 
the following inequalities:
\begin{equation}\label{ineq}
\sum_{i=1}^{j}\big(\ze_{i}(\la;\bx)-\la\big)\le \sum_{i=1}^{j}x_i\text{ for }
1\le j\le n-1.
\end{equation}
We now prove these inequalities by induction on $n$. Clearly, if 
$P(x)=x^2+ax+b$ with $a,b\in \bR$ such that $a^2>4b$ then 
$$-2\big(\ze_1(\la;\bx)-\la\big)=a+\sqrt{a^2-4b+4\la^2}\ge a+\sqrt{a^2-4b}
=-2x_1\text{ for }\la\in \bR.$$
Thus~\eqref{ineq} is true for $n=2$. Let $n\ge 3$ and assume that 
\eqref{ineq} holds for all monic strictly hyperbolic polynomials of degree at 
most $n-1$. Then we may write
\begin{equation*}
\begin{split}
& P(x)=(x-x_n)Q(x)\text{ and }P_{\la}(x)=(x-x_n)Q_{\la}(x)-\la Q(x),
\text{ where}\\
& Q(x)=\prod_{i=1}^{n-1}(x-x_i),\,Q_{\la}(x):=Q(x)-\la Q'(x)=
\prod_{i=1}^{n-1}\big(x-\om_{i}(\la;\bx')\big),\\
& \bx'=(x_1<\ldots<x_{n-1})\in \bR^{n-1}, \text{ and }\om_{i}(0;\bx')=x_i, 
\,1\le i\le n-1.
\end{split}
\end{equation*}
Note that with this labeling we get
$\om_{1}(\la;\bx')<\ldots<\om_{n-1}(\la;\bx')$ for all $\la\in\bR$ and that if 
we fix $\la$ then $\ze_{i}(\la;\bx)\rightarrow \om_{i}(\la;\bx')$ for  
$1\le i\le n-1$ while $\ze_{n}(\la;\bx)\rightarrow \infty$ as 
$x_n\rightarrow \infty$. By Lemma~\ref{ind} and the 
induction assumption applied to $Q$ we obtain
$$\sum_{i=1}^{j}\big(\ze_{i}(\la;\bx)-\la\big)< 
\sum_{i=1}^{j}\big(\om_{i}(\la;\bx')-\la\big)\le 
\sum_{i=1}^{j}x_i\text{ for }\la\neq 0\text{ and }1\le j\le n-1.$$
which proves \eqref{ineq}. We conclude that the theorem is true in this 
generic case.

{\em Step2: The general case.} Let $P$ be an arbitrary polynomial in 
$\calH_n$ and consider first an operator $T$ of the 
form $(1-\la D)e^{\la D}$ for some fixed $\la\in \bR$. As in the proof of 
Corollary~\ref{LD-span} we denote by $P_{\ve}$ the strictly 
hyperbolic polynomial with zeros $x_k+k\ve$, $1\le k\le n$, where 
$\ve>0$ and $x_1\le\ldots\le x_n$ are the (possibly multiple) zeros of $P$. 
By step 1 we have
$P_{\ve}\preccurlyeq T(P_{\ve})$ for all $\ve>0$. Using a standard
continuity argument we see that the relation $P\preccurlyeq T(P)$ is just the 
limit case $\ve\rightarrow 0$ of the aforementioned relations.
Alternatively, we may approximate $P$ with the polynomial $(1-\ve D)^{n-1}P$, 
which is strictly hyperbolic for all $\ve\neq 0$ by \cite[Lemma 4.2]{CC2}. 
Finally, if $T$ is an operator of the form $e^{-\la^2 D}$ with $\la\in \bR$ 
then 
$$P\preccurlyeq\lim_{m\rightarrow \infty}\!
\left[\left(1-\frac{\la D}{\sqrt{m}}\right)
e^{\frac{\la D}{\sqrt{m}}}\right]^{m}
\left[\left(1+\frac{\la D}{\sqrt{m}}\right)
e^{-\frac{\la D}{\sqrt{m}}}\right]^{m}\!P=T(P)$$
since $P\preccurlyeq (1-\mu D)e^{\mu D}P$ for $\mu\in \bR$. This completes 
the proof of the theorem.
\end{proof}

Note that Theorem~\ref{LP} and Corollary~\ref{LD-span} imply that operators in 
$\calA$ do not decrease the span of hyperbolic polynomials. We can actually 
deduce an even more general result from Theorem~\ref{orbit} and some 
well-known properties of classical majorization. Recall that a function 
$F:\bR^n\rightarrow \bR$ is called {\em Schur convex} if $F(X)\le F(Y)$ for 
all $X,Y\in \bR^n$ with $X\prec Y$ (cf., e.~g., \cite{MO}). Clearly, any such 
function is symmetric on $\bR^n$ and may therefore be viewed as a function on 
$\bR^n/\Sigma_n$, where $\Sigma_n$ denotes as before the symmetric group on 
$n$ elements. Theorem~\ref{orbit} and \eqref{approx} yield the 
following conditions on the relative geometry of $\calZ(P)$ and 
$\calZ\big(T(P)\big)$ for $P\in \calH_n$ and $T\in \calA'$:

\begin{corollary}\label{span-A}
Let $n\ge 2$, $P\in \calH_n$, and denote the zeros of $P$ by $x_i(P)$, 
$1\le i\le n$.
\begin{enumerate}
\item[(i)] If $T\in \calA'$ then $\min\calZ\big(T(P)\big)\le \min\calZ(P)$ and 
$\max\calZ(P)\le \max\calZ\big(T(P)\big)$, so that $\De(P)\le \De(S(P))$ for 
any operator $S\in \calA$. All these inequalities are strict unless 
$T=S=Id_{\vP}$.
\item[(ii)] The inequality $(F\circ\calZ)(P)\le (F\circ\calZ)(T(P))$ holds
for any Schur convex function $F$ on $\bR^n/\Sigma_n$ and any operator 
$T\in \calA'$. In 
particular, if $f:\bR\rightarrow \bR$ is a convex function then
$\sum_{i=1}^{n}f\big(x_i(P)\big)\le \sum_{i=1}^{n}f\big(x_i(T(P))\big)$.
\end{enumerate}
\end{corollary}

\section{Related topics and open problems}

\subsection{Toward an analytic theory of classical majorization}

Although we did not explicitly address the question of describing all 
operators in the monoid $\calA_n(\bR)$, Problem~\ref{pb1} and the results of 
the previous sections are certainly a good motivation for studying this 
question. Indeed, these results show that even a partial knowledge of 
operators in $\calA_n(\bR)$ can provide some interesting information on the 
relative geometry of the zeros of a hyperbolic polynomial and the zeros of its 
images under such operators. We therefore propose the following general 
problem.

\begin{problem}\label{pb2}
Let $n\ge 2$ and set 
$P_{\preccurlyeq}=\big\{Q\in \calH_n\mid P\preccurlyeq Q\big\}$ for 
$P\in \calH_n$. Define the monoid 
$\calB_n=\big\{T\in \calA_n(\bR)\mid P\preccurlyeq T(P)
\text{ if }P\in \calH_n\big\}$ and note that $\calA'\subseteq \calB_n$ by 
Theorem~\ref{orbit} and that $\calB_nP\subseteq P_{\preccurlyeq}$ for all 
$P\in \calH_n$, where $\calB_nP=\big\{T(P)\mid T\in \calB_n\big\}$.
\begin{enumerate}
\item[(i)] Is the inclusion 
$\calA'\subseteq \calB_n$ strict for all $n\ge 2\,$?  Describe all 
operators in $\calB_n$. 
\item[(ii)] Is it possible to describe classical majorization by means of 
the action of linear (differential) operators on hyperbolic polynomials? In 
other words, is it true that $\calB_nP=P_{\preccurlyeq}$ for all 
$P\in \calH_n$ if $n\ge 2\,$?
\item[(iii)] Characterize all operators in the monoid $\calA_n(\bR)$.
\end{enumerate}
\end{problem}

Using Corollary~\ref{span-A} it is not difficult to show that 
$\calA' P\subsetneq P_{\preccurlyeq}$ whenever $n\ge 3$ and $P\in \calH_n$ is 
such that $\De(P)>0$, where $\calA' P=\big\{T(P)\mid T\in \calA'\big\}$. 
In particular, if $P$ is strictly hyperbolic then 
$\calA' P\subsetneq P_{\preccurlyeq}$. Thus, if the answer to the first 
question in Problem~\ref{pb2} (i) were negative then it would not be possible 
to get a description of the spectral order on $\bR$ as suggested in part (ii) 
of Problem~\ref{pb2}. Nevertheless, it seems likely that 
$\calA'\subsetneq \calB_n$ for $n\ge 2$. 

Note also that an affirmative answer to Problem~\ref{pb2} (ii) would actually 
provide a description of classical majorization which in a way 
would be dual to the usual characterization by means of doubly stochastic 
matrices (cf.~Theorem~\ref{spec}): the former deals with the set 
$P_{\preccurlyeq}=\big\{Q\in \calH_n\mid P\preccurlyeq Q\big\}$ while the 
latter deals with the ``polytope'' 
$P_{\succcurlyeq}:=\big\{Q\in \calH_n\mid Q\preccurlyeq P\big\}$, where 
$P$ is an arbitrary polynomial in $\calH_n$. We refer to \cite{B2} for a 
further study of these questions and related topics.

\subsection{Pencils of complex polynomials}

The manifold $\calC_{n}$ is a 
natural context for discussing possible extensions of the results in 
\S 1--3 to the 
complex case. Indeed, by analogy with the hyperbolic case we may view 
$(\calC_{n},\preccurlyeq)$ as a partially ordered set, where the ordering 
relation $\preccurlyeq$ is now induced by the spectral order on  
$n$-tuples of vectors in $\bR^2$ (cf.~Theorem~\ref{spec} and Birkhoff's 
theorem). This means that zero sets of polynomials in $\calC_{n}$ are 
viewed as subsets of $\bR^{2}$ and that if $P,Q\in \calC_{n}$ then 
$P\preccurlyeq Q$ if and only if $Z(P)\prec Z(Q)$. The following example 
shows that if the partial ordering $\preccurlyeq$ on $\calC_{n}$ is defined 
in this way then one cannot expect a complex analog of Theorem~\ref{orbit}.

\begin{proposition}\label{counter1}
Let $P(z)=z^n-1$ and $\la \in \bC^{*}$. If $n\ge 5$ and $|\la|$ is small 
enough 
then $(1-\la D)e^{\la D}P$ and $P$ are incomparable as elements of the poset 
$(\calC_{n},\preccurlyeq)$.
\end{proposition}

\begin{proof}
Let $z_k(\la)$, $1\le k\le n$, denote the zeros of $(1-\la D)e^{\la D}P$, 
which we label so that $z_k(0)=z_k:=e^{\frac{2\pi k i}{n}}$ for
$1\le k\le n$. Since $z_k(\la)$, $1\le k\le n$, are 
analytic functions of $\la$ in a neighborhood of 0 we can use formulas
\eqref{deriv1}--\eqref{deriv2} to show
that $2\big(z_k(\la)-\la\big)=2z_k+(n-1)\bar{z}_k\la^2+\calO(\la^3)$, 
$1\le k\le n$. Hence
$$|z_k(\la)-\la|
=1+\frac{n-1}{2}\Re\big(\bar{z}_k\la^2\big)+\calO(\la^3)\text{ for } 
1\le k\le n.$$
It is geometrically clear that if $n\ge 5$ and $\la\neq 0$ then there exist 
distinct indices $k_1$ and $k_2$ such that 
$\Re\big(\bar{z}_{k_1}\la^2\big)>0$ and $\Re\big(\bar{z}_{k_2}\la^2\big)<0$. 
This implies that if $n\ge 5$ and $|\la|$ is a small enough positive number 
then $|z_{k_1}(\la)-\la|>1$ and $|z_{k_2}(\la)-\la|<1$ for some 
$1\le k_1\neq k_2\le n$. It follows that for these values of $\la$
there can be no inclusion relation between the convex hulls of 
the zeros of $(1-\la D)e^{\la D}P$ and $P$, so that these polynomials are 
incomparable as elements of the poset $(\calC_{n},\preccurlyeq)$.
\end{proof}

We also note that the results of the previous sections concerning
hyperbolic polynomials are valid only for real values of the parameter $\la$:

\begin{proposition}\label{counter2}
For any $n\ge 3$ and $\te\in (0,\pi)\cup (\pi,2\pi)$ there exists 
$\ve=\ve(n,\te)>0$ with the following property: 
if $\la=r e^{i\te}$ and $r\in (0,\ve)$ then one can find 
$P\in \calH_n$ such that $P$ and $(1-\la D)e^{\la D}P$ are 
incomparable as elements of the poset $(\calC_{n},\preccurlyeq)$.
\end{proposition}

\begin{proof}
Given a function $f:\bR^2\rightarrow \bR$ and $z\in \bC$ we shall write
$f(z)$ instead of $f(x,y)$, where $z=x+iy$. Set $I=(0,\pi)\cup (\pi,2\pi)$ 
and assume first 
that $\te\in I\setminus \{\frac{\pi}{2},\frac{3\pi}{2}\}$. Let $P$ 
be a strictly hyperbolic polynomial in $\calH_n$ with zeros $z_j$, 
$1\le j\le n$, and denote the zeros of $P-\la P'$ by $z_j(\la)$, 
$1\le j\le n$, where 
$\la=r e^{i\te}$ and $r\ge 0$. If $r$ is small enough we may label these 
zeros so that $z_j(0)=z_j$, $1\le j\le n$. Note that
$$\sum_{j=1}^{n}\frac{z_{j}P''(z_j)}{P'(z_j)}=n(n-1)\text{ and }z_j(\la)-\la=
z_j+\frac{z_{j}P''(z_j)}{2P'(z_j)}\la^2+\calO(\la^3),\,1\le j\le n,$$
by \eqref{deriv1} and \eqref{deriv2}. Fix $a\in \bR$ such that 
$|a|>|\cot 2\te|$ and consider the convex functions 
$f_{\pm a}:\bR^2\rightarrow \bR$ defined by $f_{\pm a}(x,y)=(x+ay)^2$. 
Using the above relations we get
\begin{equation*}
\begin{split}
F_{\pm a}(\la):= & \sum_{j=1}^{n}f_{\pm a}\big(z_j(\la)-\la\big)-
\sum_{j=1}^{n}f_{\pm a}(z_j)
=r^2(\cos 2\te\pm a\sin 2\te)
\sum_{j=1}^{n}\frac{z_{j}P''(z_j)}{P'(z_j)}\\
& \phantom{space}+\calO(r^3)=n(n-1)(\cos 2\te\pm a\sin 2\te)r^2+\calO(r^3).
\end{split}
\end{equation*}
Clearly, these formulas and the choice of $a$ show that there exists 
$\ve_1=\ve_1(n,\te)>0$ such that $F_{a}(\la)F_{-a}(\la)<0$ 
whenever $|\la|\in (0,\ve_1)$. By Theorem~\ref{spec} we see that
$P$ and $(1-\la D)e^{\la D}P$ are 
incomparable with respect to the partial ordering $\preccurlyeq$ on 
$\calC_{n}$.

Let us now consider the case when $\te\in \{\frac{\pi}{2},\frac{3\pi}{2}\}$, 
that is, $\la=bi$ with $b\in \bR$. Set 
$P(x)=x^n-x^{n-1}\in \calH_n$ and note that if $|\la|$ is small enough then 
there is a unique zero of $P-\la P'$ with largest real part. We denote this 
zero of $P-\la P'$
by $z(\la)$. A computation shows that $z(\la)-\la=1+(n-1)\la^2+\calO(\la^3)$ 
and so there exists $\ve_2=\ve_2(n)>0$ such that $\Re(z(\la))=1-(n-1)b^2+
\calO(b^3)<1$ if $|b|\in (0,\ve_2)$. Since $n\ge 3$ the polynomial 
$(1-\la D)e^{\la D}P$ has at least one zero at $-\la\notin \bR$. It follows 
that in this case there can be no inclusion relation between the 
convex hulls of the zeros of $(1-\la D)e^{\la D}P$ and $P$. Thus, as in the 
proof of 
Proposition~\ref{counter1} we conclude that these polynomials are 
incomparable as elements of the poset $(\calC_{n},\preccurlyeq)$. To complete 
the proof of the proposition we simply let $\ve(n,\te)=\min(\ve_1,\ve_2)$.
\end{proof}

Propositions~\ref{counter1} and~\ref{counter2} suggest that complex 
generalizations of Theorem~\ref{orbit}, if any, should involve only classical 
majorization and real values of the parameter $\la$. Using the computations 
in \eqref{approx} and \eqref{approx1} 
it is not difficult to show that if $P$ is a polynomial in $\calC_n$ whose 
zeros have distinct real parts then there exists $\ve=\ve(P)>0$ such that 
$\Re Z(P)\prec \Re Z\big((1-\la D)e^{\la D}P\big)$ if $\la$ is real and 
$|\la|\le \ve$. Based on extensive numerical experiments we make the 
following conjecture.

\begin{conjecture}\label{conj1}
If $T\in \calA'$ and $n\ge 1$ then for any $P\in \calC_n$ at least one of 
the relations $\Re Z(P)\prec \Re Z\big(T(P)\big)$, $\Im Z\big(T(P)\big)\prec
\Im Z(P)$ is valid. 
\end{conjecture}

\subsection{The Laguerre-P\'olya class of functions}

There are several known extensions of majorization to infinite sequences of 
real numbers \cite[p.~16]{MO}. A natural question is whether 
these extensions could yield infinite-dimensional analogs of 
Theorem~\ref{orbit}. For this one would need to find both a suitable
set of functions in the Laguerre-P\'olya class and an appropriate 
submonoid of $\calA'$ that acts on this set. One could for instance
consider the set of functions of genus 0 or 1 in the Laguerre-P\'olya class. 
Indeed, using Lemmas 3.1 and 3.2 in \cite{CC2} one can show that this set is 
closed under the action of operators in $\calA'$. Finally, it would be 
interesting to know whether there are any analogs of Conjecture~\ref{conj1} 
for transcendental entire functions with finitely or infinitely many 
complex zeros.

\enlargethispage{0.2cm}

\subsection{Note added in the proof}

The question of describing all linear operators that preserve the set
$\vP(\Omega)$, where $\Omega$ is a closed circular domain or the boundary
of such a domain (in particular, $\Omega=\bR$), has been settled quite recently
in \cite{BBS}.

\end{document}